\documentstyle[12pt]{article}
\author{Youssef ALAOUI}
\title{On $q$-Runge domains}
\date{}
\newcommand{\reels}{\mbox{I}\!\!\!\mbox{R}}

\newcommand{\complexes}{\mbox{I}\!\!\!\mbox{C}}

\newtheorem{th}{theorem}
\newtheorem{lm}{lemma}
\begin{document}
\maketitle
\setcounter{page}{1}
\noindent
{\large Summary}. In $[2]$, Coltoiu gave an example of a domain
$D\subset\complexes^{6}$ which is $4$-complete such that for every
${\mathcal{F}}\in Coh(\complexes^{6})$ the restriction map\\
$H^{3}(\complexes^{6},{\mathcal{F}})\rightarrow H^{3}(D,{\mathcal{F}})$
has a dense image but $D$ is not $4$-Runge in $\complexes^{6}$.
Here, we prove that for every integers $n\geq 4$ and $1\leq q\leq n$
there exists a domain $D\subset \complexes^{n}$ which is not ($\tilde{q}-1$)-Runge in $\complexes^{n}$ but such that for any
coherent analytic sheaf ${\mathcal{F}}$ on $\complexes^{n}$
the restriction map $H^{p}(\complexes^{n},{\mathcal{F}})\rightarrow H^{3}(D,{\mathcal{F}})$ has a dense image for all $p\geq \tilde{q}-2$
if $q$ does not divide $n$,
where $\tilde{q}=n-[\frac{n}{q}]+1$ and $[\frac{n}{q}]$ denotes the integral part of $\frac{n}{q}$.
\\
$2000$ MS Classification numbers: 32E10, 32E40.\\
\\
\\
{\bf 1 Introduction}\\
\\
\hspace*{.1in}An open subset $D$ of a $q$-complete space $X$
is said to be $q$-Runge in $X$ if for any compact subset  $K\subset D$
there exists a $q$-convex exhaustion function $\phi: X\rightarrow \reels$ such that $K\subset\{x\in X: \phi(x)<0\}\subset\subset D.$
When $D\subset X$ is $q$-Runge in $X$ it follows from an important result of Andreotti and Grauert $[1]$ that
for every ${\mathcal{F}}\in coh(X)$ the restriction map
$H^{p}(X, {\mathcal{F}})\rightarrow H^{p}(D, {\mathcal{F}})$ has a dense image for all $p\geq q-1$.\\
\hspace*{.1in}A counterexample to the converse of this statement is given by Coltoiu in $[2]$ where it is shown that there exists a domain
$D\subset \complexes^{6}$ which is $4$-complete but not $4$-Runge in $\complexes^{6}$ such that for every
${\mathcal{F}}\in Coh(\complexes^{6})$ the restriction map
$H^{3}(\complexes^{6},{\mathcal{F}})\rightarrow H^{3}(D,{\mathcal{F}})$
has a dense image.
\newpage
\noindent
\hspace*{.1in}The purpose of this paper is to prove the following.
\begin{th}{-Let $(n,q)$ be a pair of integers with $1\leq q\leq n$.
Suppose that $n\geq 4$ and $q$ does not divide $n$. Then there is a domain $D\subset \complexes^{n}$ which is not $(\tilde{q}-1)$-Runge in $\complexes^{n}$ but for any coherent analytic sheaf ${\mathcal{F}}$ on $\complexes^{n}$ the restriction map
$H^{p}(\complexes^{n},{\mathcal{F}})\rightarrow H^{p}(D,{\mathcal{F}})$
has a dense image for all $p\geq \tilde{q}-2$.}
\end{th}
{\bf Examples proving theorem 1}\\
\\
\hspace*{.1in}Let $(n,q)$ be a pair of integers with $1\leq q\leq n.$
We put $m=[\frac{n}{q}]$ and suppose that $r=n-mq>0.$ We consider
the functions $\phi_{1},\cdots, \phi_{m+1}$ defined on $\complexes^{n}$
by
$$\phi_{j}=\sigma_{j}+\displaystyle\sum_{i=1}^{m}\sigma_{i}^{2}
-\frac{1}{4}||z||^{2}+N||z||^{4}, \ \ j=1,\cdots,m,$$
and
$$\phi_{m+1}=-\sigma_{1}-\cdots-\sigma_{m}+\displaystyle\sum_{i=1}^{m}\sigma_{i}^{2}
-\frac{1}{4}||z||^{2}+N||z||^{4},$$
where $\sigma_{j}=y_{j}+\displaystyle\sum_{i=m+1}^{n}|z_{j}|^{2}-
(m+1)\displaystyle\sum_{i=m+(j-1)(q-1)+1}^{m+j(q-1)}|z_{j}|^{2}$
for $j=1,\cdots,m$, $z_{j}=x_{j}+iy_{j}$. Then it is known from $[4]$
that all $\phi_{j}$, $1\leq j\leq m+1$, are $q$-convex on $\complexes^{n}$, if $N>0$ is sufficiently large and, if
$\rho=Max\{\phi_{j}, 1\leq j\leq m+1\}$, then for $\varepsilon_{o}>0,$
the set
$D_{\varepsilon_{o}}=\{z\in\complexes^{n}: \rho(z)<-\varepsilon_{o}\}$
is relatively compact in the unit ball $B=B(0,1)$, if $N$ is large enough and $\varepsilon_{o}>0$ is sufficiently small.
\begin{lm}{-For every coherent analytic sheaf ${\mathcal{F}}$ on
$\complexes^{n}$, the cohomology groups
$H^{p}(D_{\varepsilon_{o}},{\mathcal{F}})$ vanish for all $p\geq
\tilde{q}-2$.}
\end{lm}
{\bf Proof}
\\
\hspace*{.1in}The proof is by induction on $m$.\\
\hspace*{.1in}We have $D_{\varepsilon_{o}}=D_{1}\cap\cdots\cap D_{m+1}$,
where $D_{i}=\{z\in \complexes^{n}: \phi_{i}(z)<-\varepsilon_{o}\}$
is obviously $q$-complete and $q$-Runge in $\complexes^{n}$. In fact,
$\frac{-1}{\phi_{i}+\varepsilon_{o}}$ is a $q$-convex exhaustion function on $D_{i}$
and, if $K\subset D_{i}$ is a compact subset and\\ $\psi\in C^{\infty}(\complexes^{n},(0,\infty))$ a strictly plurisubharmonic exhaustion function on $\complexes^{n}$, then for every real number $C>0,$ the function $\psi_{C}=\psi+C(\phi_{i}+\varepsilon_{o})$ is a $q$-convex  exhaustion on $\complexes^{n}$ and, if $C$ is sufficiently large, then
$$K\subset\{z\in\complexes^{n}: \psi_{K}(z)<0\}\subset\subset D_{i}$$
\hspace*{.1in}If now $m=1$, then $\tilde{q}=n$ and $D_{\varepsilon_{o}}=D_{1}\cap D_{2}$,
where each $D_{i}$ is, in particular, $(n-1)$-Runge in $\complexes^{n}$. Then the restriction map $H^{n-2}(\complexes^{n},{\mathcal{F}})\rightarrow H^{n-2}(D_{i},{\mathcal{F}})$ has a dense image. This proves
that $H^{n-2}(D_{i},{\mathcal{F}})=0$ for $i=1, 2$. Because $\complexes^{n}\backslash D_{1}\cup D_{2}$ has no compact connected components, then, by [3],
$D_{1}\cup D_{2}$ is $n$-Runge in $\complexes^{n}$, which implies
that $H^{n-1}(D_{1}\cup D_{2}, {\mathcal{F}})=0$. Therefore
from the Mayer-Vietoris sequence for cohomology
\begin{center}
$\cdots\rightarrow H^{n-2}(D_{1},{\mathcal{F}})\oplus
H^{n-2}(D_{2},{\mathcal{F}})\rightarrow H^{n-2}(D_{\varepsilon_{o}},{\mathcal{F}})\rightarrow H^{n-1}(D_{1}\cup D_{2},{\mathcal{F}})\rightarrow\cdots$
\end{center}
it follows that $H^{n-2}(D_{\varepsilon_{o}},{\mathcal{F}})=0.$\\
\hspace*{.1in}Suppose now that for any integer $k$ with $1\leq k\leq m$
the cohomology group $H^{p}(D_{i_{1}}\cap\cdots\cap D_{i_{k}},{\mathcal{F}})=0$ for all
$p\geq \tilde{q}-2$ and $i_{1},\cdots, i_{k}\in\{1,2,\cdots,m+1\}$.
Then, by [5, Proposition $1$],
$H^{\tilde{q}-2}(D_{1}\cap\cdots\cap D_{m+1},{\mathcal{F}})\simeq
H^{n-1}(D_{1}\cup\cdots\cup D_{m+1},{\mathcal{F}}).$
Since $\complexes^{n}\backslash D_{1}\cup\cdots\cup D_{m+1}$
has no compact connected components, then $D_{1}\cup\cdots\cup D_{m+1}$
is $n$-Runge. Therefore
$H^{\tilde{q}-2}(D_{1}\cap\cdots\cap D_{m+1},{\mathcal{F}})=
H^{n-1}(D_{1}\cup\cdots\cup D_{m+1},{\mathcal{F}})=0$.
\begin{lm}{-The set $D_{\varepsilon_{o}}$ is not $(\tilde{q}-1)$-Runge in $\complexes^{n}$}
\end{lm}
{\bf Proof}\\
\\
\hspace*{.1in}There exists small enough $\delta>0$ such that the following topological sphere of real dimension $2n-m-1$:
\begin{center}
$S_{\delta}=\{z\in\complexes^{n}:
x_{1}^{2}+\cdots+x_{m}^{2}+|z_{m+1}|^{2}+\cdots+|z_{n}|^{2}=\delta,
y_{j}=-\displaystyle\sum_{i=m+1}^{n}+(m+1)\displaystyle\sum_{i=m+(j-1)(q-1)+1}^{m+j(q-1)}|z_{i}|^{2} \ for \ j=1,\cdots, m\}$
\end{center}
is contained in $D_{\varepsilon_{o}}$. It was shown by Diederich-Fornaess $[4]$ that $S_{\delta}$ is not homologous to $0$ in $D_{\varepsilon_{o}}$, because $D_{\varepsilon_{o}}$ does not meet\\
$E=\{z\in \complexes^{n}: x_{1}=\cdots=x_{m}=z_{m+1}=\cdots=z_{n}=0\}.$\\
\hspace*{.1in}Suppose now that $D_{\varepsilon_{o}}$ is $(\tilde{q}-1)$-Runge in $\complexes^{n}$. Then there exists a smooth
$(\tilde{q}-1)$-convex exhaustion function $\phi$ on $\complexes^{n}$
such that
$$S_{\delta}\subset \tilde{D}_{o}=\{z\in\complexes^{n}: \phi(z)<0\}\subset\subset D_{\varepsilon_{o}}.$$
Since $\tilde{D}_{o}\cap E=\emptyset$, then $S_{\delta}$ is not homologous to $0$ in $\tilde{D}_{o}$. Moreover,
there exists a sufficiently large constant $C>0$ such that
$$B\subset \tilde{D}_{C}=\{z\in \complexes^{n}: \phi(z)<C\}$$
But since in addition $\phi$ is $(\tilde{q}-1)$-convex, it follows
from $[6]$ that
$$H_{2n-(m+1)}(\tilde{D}_{o},\reels)\simeq H_{2n-(m+1)}(\tilde{D}_{C},\reels),$$
which is false, because $S_{\delta}$ is clearly homologous to $0$
in $\tilde{D}_{C}$.

\end{document}